\documentclass[12pt]{article}
\usepackage{enumerate}
\usepackage{amsmath}
\usepackage{amsthm}
\usepackage{amssymb}
\usepackage{color}
\usepackage{graphicx}
\usepackage{xspace}
\usepackage{url}

\newtheorem{theorem}{Theorem}

\newcommand{\Ppresent}[2]{\noindent\textbf{Presented by #1, 
\texttt{#2}\\}}
\newcommand{\Porig}[1]{\noindent\textbf{Original proposer: #1.\\}}
\newcommand{\av}{\mathop{\mathrm{av}}\nolimits}

\begin{document}

\title{Problems from BCC30}
\author{Edited by Peter Cameron, St Andrews and QMUL}
\date{\texttt{pjc20@st-andrews.ac.uk}}
\maketitle
\begin{abstract}
These problems were mostly presented at the problem session at the 30th
British Combinatorial Conference at Queen Mary University of London on
4 July 2024. Some were contributed later by conference participants. Thank
you to all the contributors.

The problems are given here in alphabetical order of presenter.
If no originator is given, I assume that the presenter is the originator.

Please send corrections to me (\texttt{pjc20@st-andrews.ac.uk}). Solutions
should be sent to the presenter; I would appreciate a copy too.
\end{abstract}

\setcounter{section}{30}

\subsection{Inversions and avoidance}
\Ppresent{Christian Bean}{c.n.bean@keele.ac.uk}
\Porig{Anders Claesson, V\'{\i}t Jel\'{\i}nek and Einar Stein\-gr\'{\i}msson}

An \emph{inversion} in a permutation $\pi$ is a pair $(i,j)$ of indices
such that $i<j$ and $\pi(i)>\pi(j)$.

Claesson, Jel\'{\i}nek and Steingr\'{\i}msson~\cite{cjs} conjectured that, in
the class $\av(1324)$ of permutations avoiding the pattern $1324$, for
fixed~$k$ the number of permutations with exactly $k$ inversions grows
monotonically with the length of the permutation. They also showed 
that the truth of the conjecture would imply that the growth constant for
this class (the limit of $\av_n(1234)^{1/n}$ as $n\to\infty$) is
$\mathrm{e}^{\pi\sqrt{2/3}}=13.002\ldots$.

\paragraph{Problem} Prove the conjecture.

\subsection{Making a tournament transitive}
\Ppresent{Natalie Behague}{natalie.behague@warwick.ac.uk}
\newcommand{\inv}{\mathop{\mathrm{inv}}}

For a tournament $T$, let $\inv(T)$ be the minimum number of subsets
$X_1,\ldots,X_k$ of the vertex set $V(T)$ with the property that inverting
the edges within $X_i$ for all $i$ in turn gives a \emph{transitive}
tournament (a total order)?

Bang-Jensen, Costa Ferreira da Silva and Havet~\cite{bch} conjectured that,
for all $T_1$ and $T_2$,
\[\inv(T_1\to T_2)=\inv(T_1)+\inv(T_2),\]
where $T_1\to T_2$ is the disjoint union of $T_1$ and $T_2$ together with
all edges directed from $T_1$ to $T_2$. However, this conjecture was refuted
by Alon, Powierski, Savery, Scott and Wilmer~\cite{apssw}, and by
Aubian, Havet, H\"orsch, Klingelhoefer, Nisse, Rambaud and Vermande
\cite{ahhknrv}.

\paragraph{Question} Is there a constant $c$ such that
\[\inv(T_1\to T_2) \ge \inv(T_1) + \inv(T_2)-c?\]
Perhaps $c=1$ or $c=2$ \dots

\medskip

There is a link to matrix rank. See Behague, Johnston, Morrison and
Ogden~\cite{bjmo}, and Wang, Yang and Lu~\cite{wyl}.

\subsection{A game on the edges of $K_n$}
\Ppresent{Thomas Bloom}{bloom@maths.ox.ac.uk}
\Porig{Paul Erd\H{o}s}

This Erd\H{o}s problem is taken from~\cite{guy}.

This is a game played on the edge set of the complete graph $K_n$. There
are two players, {\color{red}RED} and {\color{blue}BLUE}. {\color{red}RED}
goes first. The players colour the edges alternately until no uncoloured
edges remain. {\color{red}RED} wins if the largest red clique is strictly
larger than the largest blue clique.

\paragraph{\color{black}Conjecture} {\color{blue}BLUE} wins, for all $n\ge3$.

\subsection{A conjectured new version of the Two Families Theorem}
\Ppresent{Asier Calbet}{a.calbetripodas@qmul.ac.uk}

The Two Families Theorem is a celebrated result in extremal set theory. There
are several different versions of this theorem, the simplest one being as
follows.

\paragraph{Theorem} (Two Families Theorem, Bollob\'as, 1965).
Let $a, b\ge 0$ be integers and $(A_i , B_i )_{i\in I}$ be a
collection of pairs of finite sets, indexed by a finite set $I$, with the
following properties:
\begin{itemize}
\item $|A_i| = a$, $|B_i| = b$ and $A_i\cap B_i = \emptyset$ for all $i\in I$;
\item $A_i\cap B_j\ne\emptyset$ for all $i,j\in I$ with $i\ne j$.
\end{itemize}
Then
\[|I|\le{a+b\choose a}.\]

\medskip

One can see that Theorem 1 is tight by taking $(A_i,B_i)_{i\in I}$ to be the
collection of all partitions of a set of size $a+b$ into a subset $A_i$ of
size $a$ and a subset $B_i$ of size $b$. Moreover, this is the unique way of
achieving equality. There are numerous extensions and variations of this
theorem and it has many applications.

In~\cite{calbet}, I prove a new version of the Two Families Theorem which I
then use to prove bounds on the minimum number of edges in a $K_r$-saturated
graph with given minimum degree. I believe that to prove tight bounds for this
saturaton problem a further generalisation of the Two Families Theorem is
needed. The simplest special case of this generalisation that I cannot
prove is as follows.

\paragraph{Conjecture}
Let $b\ge a\ge 2$ be integers and $(A_i,B_i)_{i\in I}$ be a collection of
pairs of finite sets, indexed by a finite set $I$, with the following
properties:
\begin{itemize}
\item $|A_i|=a$, $|B_i|=b$, and $|A_i\cap B_i|=2$ for all $i\in I$;
\item $A_i\cap B_j\not\subseteq A_k\cap B_k$ for all $i,j,k\in I$
with $i\ne j$.
\end{itemize}
Then
\[|I|\le\sum_{i=2}^a2^{i-2}{a+b-2i\choose a-i}.\]

\medskip

The following example shows that, if true, Conjecture 1 is tight. For each
integer $c$ with $2\le c\le a$, let $(A_i)_{i\in I_c}$ be the
collection of all subsets $A_i$ of $\{1,\ldots,a+b-2\}$ satisfying
\begin{itemize}
\item $|A_i|=a$;
\item $2c-3,2c-2\in A_i$;
\item $|A_i\cap\{2d-3,2d-2\}|=1$ for all $d$ with $2\le d<c$.
\end{itemize}
For each $i\in I_c$, let $B_i=(\{1,\ldots,a+b-2\}\setminus A_i)\cup
\{2c-3,2c-2\}$. Let $I=\bigcup_{2\le c\le a}I_c$. Then
$(A_i,B_i)_{i\in I}$ satisfies the conditions in Conjecture~1 and meets
the conjectured bound.

I can prove Conjecture~1 under the additional assumption that the sets
$A_i$ and $B_i$ are all subsets of a ground set of size $a+b-2$ (the
smallest possible size of a set containing sets of sizes $a$ and $b$ whose
intersection has size~$2$). Note that taking $k=i$ in the last condition,
we obtan that $(A_i\setminus B_i,B_i)_{i\in I}$ satisfies the conditions
in Theorem 1 with $a$ replaced by $a-2$. Hence $|I|\le{a+b-2\choose a-2}$.
When $a=3$, this bound matches the one in Conjefcture~1 but there are many
different ways of achieving equality.

See~\cite{calbet} for more details. In particular, Problem~2 and 
Definition~11 of that paper give the full generalisation of the Two
Families Theorem needed for the saturation application. See also
\cite[Section 4]{scott_wilmer} for a discusson of several different
versions of the Two Familires Theorem and their applications.

\subsection{Zero link Tur\'an densities}
\Ppresent{Asier Calbet}{a.calbetripodas@qmul.ac.uk}
\newcommand{\ex}{\mathop{\mathrm{ex}}}

This problem arose after reading the final remark in the recent paper
\cite{eil}.

For an $r$-graph $H$, let
\[\pi_r(H)=\lim_{n\to\infty}\ex(n,H)/{n\choose r}\]
be its Tur\'an density, which exists by the usual averaging argument. We
define the \emph{link Tur\'an density} of $H$ based on~\cite{eil}. For
every integer $k\ge0$, let $H_k$ be the $(r+k)$-graph obtained from $H$
by adding $k$ vertices and replacing each edge by its union with this set
of $k$ vertices. Note that an $(r+k)$-graph contains $H_k$ if and only if
one of its $k$-links contains $H$. Then we define the link Tur\'an density as
\[\pi_\infty(H)=\lim_{k\to\infty}\pi_{r+k}(H_k),\]
which exists by a similar averaging argument.

The obvious problem is to determine $\pi_{\infty}(H)$, but this may be very
difficult. So an easier, but still interesting, question to ask is 

\paragraph{Problem} For which $r$-graphs $H$ do we have $\pi_\infty(H)=0$?

\medskip

It is known that, for an $r$-graph $H$, $\pi_r(H)=0$ if and only if $H$ is
$r$-partite. Perhaps the above problem has a similar answer.

\medskip

The final remark in \cite{eil} is precisely the solution to the Problem when
$r=2$:

\paragraph{Theorem} For a $2$-graph $H$, $\pi_\infty(H)=0$ if and only if $H$
is tripartite.

\medskip

The case $r=3$ is still open. Note that, since $\pi_\infty(H)=\pi_\infty(H_1)$,
the problem becomes harder as $r$ increases.

The \emph{dual} $H^*$ of an $r$-graph $H$ on $v$ vertices is the $(v-r)$-graph
on the same vertex set whose edges are the complements of the edges of $H$.
Now it is possible to show:

\paragraph{Lemma} The set of hypergraphs $H$ with $\pi_\infty(H)=0$ is closed
under subgraphs, blow-ups and duals.

\subsection{Cliques and cocliques in the Johnson scheme}
\Ppresent{Peter Cameron}{pjc20@st-andrews.ac.uk}
\Porig{Mohammed Aljohani, John Bamberg, Peter J. Cameron}

Let $I$ be a non-empty proper subset of $\{0,\ldots,k-1\}$. Let $A$ and $B$
be sets of $k$-subsets of $\{1,\ldots,n\}$. Assume that
\clearpage
\begin{itemize}
\item for $a_1,a_2\in A$, we have $|a_1\cap a_2|\in I$;
\item for $b_1,b_2\in B$, we have $|b_1\cap b_2|\in\{0,\ldots,k-1\}\setminus I$;
\item $|A|\cdot|B|={n\choose k}$.
\end{itemize}
(The first two bullet points imply that $|A|\cdot|B|\le{n\choose k}$.)

\paragraph{Problem} Show that, if $n$ is sufficiently large in terms of $k$,
then there exists $t$ such that either $I$ or $\{0,\ldots,k-1\}\setminus I$ is
equal to $\{0,1,\ldots,t-1\}$.

\medskip

This is true for $k\le4$, see \cite{abc}.

Note that, if $I=\{0,\ldots,t-1\}$, then $A$ is the set of blocks of a Steiner
system $S(t,k,n)$, and $B$ is a maximum-size $t$-intersecting family.

\subsection{Compatible orderings of vector spaces}
\Ppresent{Peter Cameron}{pjc20@st-andrews.ac.uk}

A total ordering of a vector space is \emph{compatible} if the unique
order-preserving bijection between any two subspaces of the same dimension
is linear.

It is not hard to show that every finite vector space has a compatible
ordering. For example, choose any total ordering of the field, and then
order the vector space lexicographically.

\paragraph{Problem} How many compatible orderings does the vector space of
dimension $n$ over the $q$-element field have?

\medskip

The answer is known for $n=1$ and for $n=2$.

\subsection{Random walks on graphs on groups}
\Ppresent{Peter Cameron}{pjc20@st-andrews.ac.uk}

Let $G$ be a finite group. The \emph{commuting graph} of $G$ is the graph
with vertex set $G$, in which $x$ and $y$ are joined if $xy=yx$ (this
includes a loop at each vertex). It may have been Mark Jerrum who first
pointed out that the limiting distribution of the random walk on the
commuting graph is \emph{uniform on conjugacy classes}, that is, the
probability of an element $g$ is inversely proportional to the size of the
conjugacy class of $g$.

\paragraph{Problem}
Describe the limiting distributions of the random walk on other graphs on
groups. Examples of interest include
\begin{itemize}
\item the \emph{generating graph} ($x$ and $y$ are joined if the pair $x,y$
generates $G$); or
\item the \emph{power graph} ($x$ and $y$ are joined if one is a power of the
other).
\end{itemize}

\subsection{Associated Stirling numbers}
\Ppresent{Bishal Deb}{bishal@gonitsora.com}

Let $\genfrac[]{0pt}0{n}{k}_{\ge r}$ be the number of permutations in the
symmetric group $\mathfrak{S}_n$ which have $k$ cycles, each of length at
least~$r$; these are known as the \emph{$r$-associated Stirling cycle numbers}
\cite[pp.256--257,295]{comtet}. The \emph{$r$-th order Stirling cycle
numbers}~\cite[Ch.~6]{deb} are given by
\[\genfrac[]{0pt}0{n}{k}^{(r)}=\genfrac[]{0pt}0{n+(r-1)k}{k}_{\ge r}.\]
Define the \emph{$r$-th order Stirling cycle polynomials} $c_{r,n}(x)$ to be
\[c_{r,n}(x)=\genfrac[]{0pt}0{n}{k}^{(r)}x^k.\]

\paragraph{Fact}~\cite[Theorem 6.1.6]{deb}
For $r=1,2$, the polynomials $c_{r,n}(x)$ are real-rooted,
and so the sequence
\[\genfrac[]{0pt}0{n}{0}^{(r)},\genfrac[]{0pt}0{n}{1}^{(r)},\ldots,
\genfrac[]{0pt}0{n}{n}^{(r)}\]
is log-concave.

\paragraph{Conjecture}~\cite[Conjecture 6.1.2(c)]{deb}
It is also log-concave for $r=3,4,5$.

\clearpage

\subsection{Permutation patterns with squares}
\Ppresent{Stoyan Dimitrov}{emailtostoyan@gmail.com}

\paragraph{Definition}
A permutation $\pi=\pi_1\pi_2\cdots\pi_n\in S_n$ \emph{contains} the sequence
$\sigma=\sigma_1\sigma_2\cdots\sigma_t\Box\sigma_{t+1}\cdots\sigma_k$
\emph{as a pattern} if
\begin{itemize}
\item there exist $x_1,\ldots,x_k$ with $1\le x_1<\cdots<x_k\le n$ such that,
for all $i,j\in\{1,\ldots,k\}$ with $i\ne j$, we have $\pi_{x_i}<\pi_{x_j}$
if and only if $\sigma_i<\sigma_j$;
\item $x_{t+1}-x_t>1$ (that is, the $\Box$ symbol indicates a gap of size
at least~one).
\end{itemize}
If $\pi$ does not contain $\sigma$, we say that $\pi$ \emph{avoids} $\sigma$.

\medskip

For example, $\pi=52341$ contains $\sigma=31\Box2$, because of its subsequence
$\pi_1\pi_2\pi_4=524$; but it does not contain the pattern $32\Box1$, because
$541$ is the only occurrence of $321$, but the letters $4$ and $1$ are
adjacent. Denote by $\av_n(p)$ the number of permutations in $S_n$
avoiding $p$. Two patterns $p$ and $q$ (possibly with squares) are
\emph{Wilf-equivalent} if $\av_n(p)=\av_n(q)$ for every $n\ge1$.

\paragraph{Conjecture} \cite[Conjecture 8.1]{dimitrov} Take the Wilf-equivalent
patterns  $1234$, $1243$ and $2143$ and put a square in the same position in
each of these patterns; the results will be Wilf-equivalent. That is, for every
$n\ge1$, we have
\begin{eqnarray*}
&&\av_n(1\Box234)=\av_n(1\Box243)=\av_n(2\Box143),\\
&&\av_n(12\Box34)=\av_n(12\Box43)=\av_n(21\Box43),\\
&&\av_n(123\Box4)=\av_n(124\Box3)=\av_n(214\Box3).
\end{eqnarray*}

\medskip

\paragraph{Note} A similar statement does not hold for every two
Wlf-equivalent patterns, because $3142$ and $4132$ are Wilf-equivalent
\cite{stankova}, but
\[av_7(4\Box132)=3592\ne3587=\av_7(3\Box142).\]

\clearpage

\subsection{From a rotation system to a surface}
\Ppresent{Andrei Gagarin}{GagarinA@cardiff.ac.uk}
\Porig{Andrei Gagarin, William Kocay}

We refer to \cite{gagarin_kocay} for definitions.

Given a rotation system of a $2$-cell embedding of a $3$-connected graph on a
topological surface, either orientable or non-orientable, find a corresponding
drawing of the graph on the topological surface.

This problem has been solved for the plane, torus, and projective plane, as
described in the book~\cite{kk}. Software is available for finding drawings
on these three surfaces.

The problem is currently unsolved for surfaces of higher genus.
For example, as mentioned in \cite[p.~41]{gagarin_kocay}, there are $13$
rotation systems for $K_5$ on the triple torus. I have managed to draw two
of them, shown in Figure 23 of the paper. So, a starting point would be to
draw the other $11$ embeddings given by the rotation systems on the triple
torus. Note that the covering surfaces of the torus and projective plane are
the Euclidean plane and the 2-sphere. Yet beginning with the double and triple
tori, the covering surface is the hyperbolic plane.

\paragraph{Problem} Draw all the $13$ different rotation systems of $K_5$ on
the triple torus.

\subsection{A geolocation problem}
\Ppresent{\'Alvaro Guti\'errez Caceres}{\hfil\break 
a.gutierrezcaceres@bristol.ac.uk}

Let $\mathcal{C}$ be a compact subspace of $\mathbb{R}^n$. A point
$x\in\mathcal{C}$ is determined uniquely by its distances to $n+1$ points
in general position. The same holds if $x$ is replaced by a sphere $X$
(Appolonus' Theorem). We can turn this into a one-player game, where the
player selects a point $a_k\in\mathbb{R}^n$ at each turn, and an oracle
returns the distance $f(a_1,\ldots,a_k)=d(a_k,X)$ to the sphere.

\paragraph{Problem}
Changing the oracle so that it returns
\[f(a_1,\ldots,a_k)=\min_{i\le k} d(a_i,X),\]
how many turns are required to determine $X$? Give an algorithm to select
$a_1,a_2,\ldots$.

\medskip
Proven lower bounds for the number of turns: $n+1$.
Proven upper bounds for the number of turns:
\begin{itemize}
\item $3$ if $n=2$;
\item $n^2+n+1$ for general $n$.
\end{itemize}

\paragraph{Conjecture} The lower bound is sharp.

Known results can be found in \cite{gutierrez}.

\subsection{Supertrees}
\Ppresent{Zach Hunter}{zachtalkmath@gmail.com}

We say that a tree $T$ is an \textit{$n$-supertree} if any $n$-vertex tree
$T'$ can be obtained from $T$ by edge contractions (so that $T'$ is a minor
of $T$).

Let $F(n)$ be the smallest number of vertices of an $n$-supertree.

\paragraph{Problem} Is $F(n)\le n^{1+o(1)}$?

\medskip

It is known that $cn\log n\le F(n)\le n^2$. It suffices to assume that
$\Delta(T')\le3$.

\subsection{Permutations shattering many triples}
\Ppresent{Robert Johnson}{r.johnson@qmul.ac.uk}
\Porig{J.~Robert Johnson and Belinda Wickes}

This is an extremal question about permutations. For more details see
\cite{johnson_wickes}.

A family $\mathcal{P}$ of permutations from $S_n$ {\em shatters} a $k$-set $X \subset [n]$ if the elements of $X$ appear in each of the
$k!$ possible orders in permutations from $\mathcal{P}$.

For example, the following family of permutations from $S_5$ shatters the triple $\{2, 3, 5\}$, but not $\{1, 2, 3\}$, since no permutation gives the order $(1, 3, 2)$.
\[
\{12345, 35241, 41523, 25143, 53142, 43215\}
\]

In fact this family shatters $8$ of the $10$ triples from $[5]$ and it can be checked (by a rather tedious case analysis) that this is the largest possible; there is no family of $6$ permutations of $[5]$ which shatters $9$ triples.

\paragraph{Problem}
What is the largest number of triples from $[n]$ which can be shattered by a family of $6$ permutations?

\medskip

Writing $f(n)$ for this quantity, it is easy to see by an averaging argument that $\frac{f(n)}{\binom{n}{3}}$ is decreasing in $n$ and so this ratio tends to a limit. We have the following bounds:
\[
\frac{17}{42}\leq\lim_{n\to\infty}\frac{f(n)}{\binom{n}{3}}\leq\frac{47}{60}.
\]

\subsection{Package colouring}
\Ppresent{George Kontogeorgiou}{gkontogeorgiou@dim.uchile.cl}
\Porig{Martin Winter}

A map $c:V(G)\to\mathbb{N}$ on the vertex set of a graph $G$, such that
$c(V(G))$ is finite, is said to be a \emph{package colouring} if any two
distinct vertices $u,v$ with $c(u)=c(v)=n$ satisfy $d(u,v)\ge n+1$.

Here is an example.
\begin{center}
\setlength{\unitlength}{1mm}
\begin{picture}(90,10)
\put(0,0){\line(1,0){90}}
\multiput(0,0)(10,0){8}{\circle*{1}}
\multiput(-0.5,1.5)(20,0){4}{$1$}
\multiput(9.5,1.5)(40,0){2}{$2$}
\multiput(29.5,1.5)(40,0){2}{$3$}
\end{picture}
\end{center}

For a subset $S$ of $\mathbb{N}$, put
\[\partial(S)=\sum_{s\in S}\frac{1}{s+1}.\]
Then, for example, $\partial(\{2^i:i\ge0\})=1.26\ldots$.

If $\partial(S)\ge2$ and $S$ is finite, then $S$ package colours $\mathbb{Z}$.

\paragraph{Problem} What is
\[\inf\{x\in\mathbb{R}\mid\partial(S)\ge x\Rightarrow S\hbox{ package colours }
\mathbb{Z}\}?\]

\clearpage

\subsection{Cuboctahedra in Latin squares}
\Ppresent{Matthew Kwan}{matthew.kwan@ist.ac.at}

\paragraph{Definition} A \emph{cuboctahedron} in a Latin square $L$ is a
pair of pairs of rows $(r_1,r_2)$ and $(r_1',r_2')$ and a pair of pairs of
columns $(c_1,c_2)$ and $(c_1',c_2')$ such that $L(r_i,c_j)=L(r_i',c_j')$
for all $i,j\in\{1,2\}$.

\begin{theorem}[Brandt~\cite{brandt}] $L$ is the Cayley table of a group if
and only if it has $n^5$ cuboctahedra (the maximum number possible).
\end{theorem}

\paragraph{Editor's note} This is related to the theorem of
Frolov~\cite{frolov}.

\paragraph{Observation} Every Latin square has at least $(3-o(1))n^4$
cuboctahedra.

\paragraph{Question} Is there a Latin square with this few cuboctahedra?

\begin{theorem}[Kwan--Sah--Sawhney--Simkin~\cite{ksss}]
A random Latin square has $(4+o(1))n^4$ cuboctahedra.
\end{theorem}

\subsection{Flip graphs}
\Ppresent{Xandru Mifsud}{xandru.mifsud.16@um.edu.mt}
\Porig{Yair Caro, Josef Lauri, Xandru Mifsud, Raphael Yuster and Christina
Zarb}

For a graph $G = (V, E)$ which is edge-coloured using colours $[k] = \{1, \dots, k\}$, given a vertex $v \in V$ and $j \in [k]$, by $e_j [v]$ we denote the number of edges coloured $j$ in the closed neighbourhood of $v$, whilst $\deg_j (v)$ denotes the number of edges coloured $j$ incident to $v$.

Consider an increasing positive-integer sequence $(a_1, \dots, a_k)$ and let $d = \sum_{j=1}^k a_j$. A $d$-regular graph $G = (V, E)$ is said to be an \textit{$(a_1,\ldots,a_k)$-flip graph} if there exists an edge-colouring $ f : E \rightarrow  \{1, \dots, k\}$ such that
\begin{itemize}
\item the edges coloured $j$ span an $a_j$-regular subgraph, namely
$\deg_j(v) = a_j$ for every $v\in V$, resulting in a global majority ordering;
\item for every vertex $v \in  V$,  $e_k[v] < e_{k-1}[v] < \ldots <  e_1[v]$,
resulting in a locally opposite majority ordering.
\end{itemize} 
	
We say that $(a_1,\ldots,a_k)$ is a \textit{flip sequence} if it is realised by some $(a_1,\ldots,a_k)$-flip graph. It is known that $(a_1, a_2)$ is a flip sequence if, and only if, $3 \leq a_1$ and $a_2 < \binom{a_1 + 1}{2}$. Moreover if $(a_1, a_2, a_3)$ is a flip sequence, then $a_3 < 2(a_1)^2$. 

However for flip sequences on 4 or more colours, the situation changes drastically: Let $k \in \mathbb{N}$ such that $k > 3$. Then there is some constant $m = m(k) \in \mathbb{N}$ such that for all $N \in \mathbb{N}$, there exists a $k$-flip sequence $\left(a_1, a_2, \dots, a_k\right)$ such that $a_1 = m$ and $a_k > N$. 

In other words, the largest colour degree is now no longer quadratically bound in terms of the smallest!

\paragraph{Problem}
	Let $k \in \mathbb{N}$ such that $k > 3$. What is the largest $q(k) \in \mathbb{N}$, $q(k) < k$, for which there exists some positive-integers $a_1 < a_2 < \dots < a_{q(k)}$ such that for all $N \in \mathbb{N}$ there exists an $(a_1, \dots, a_{q(k)}, \dots a_k)$-flip sequence of length $k$ satisfying $a_k > N$?

\medskip

It is known that for $k > 3$ we have that $$\max\left\{1, \left\lceil\frac{k}{4}\right\rceil - 1\right\} \leq q(k) < \begin{cases}
 	\ \frac{k}{3} & \mbox{if $k \equiv 0 \ (\!\!\!\!\!\!\mod 3)$} \\
 	\left\lceil \frac{k}{2}\right\rceil & \mbox{otherwise}
 \end{cases}
$$ where the lower-bound is sharp. 

\subsection{Graphs containing every tree}
\Ppresent{Nika Salia}{salianika@gmail.com}

Let $S(n)$ be the minimum number of edges in a graph on $n$ vertices which
contains every tree on $n$ vertices as a subgraph.

It is known that 
\[n\log n-O(n) \le S(n) \le \frac{5}{\log 4}n\log n+O(n).\]
The upper bound is due to Chung and Graham~\cite{chung_graham}, the lower
bound to Gy\H{o}ri, Li, Salia and Tompkins~\cite{glst}.

\paragraph{Problem} Find $S(n)$.

\subsection{Cycles in bipartite graphs}
\Ppresent{Nika Salia}{salianika@gmail.com}

Let $G$ be a bipartite graph with bipartite parts $A$ and $B$.

For a subset $D$ of $A$, let $\hat N(D)$ be the set of vertices having at
least two neighbours in $D$.

\paragraph{Conjecture} Suppose that for all $D\subset A$ with $|D|\ge2$,
we have $|\hat N(D)|\ge|D|$. Then $G$ contains a cycle of length $2|A|$.

\subsection{An orientation question}
\Ppresent{Thomas Selig}{Thomas.Selig@xjtlu.edu.cn}

Let $G=(V,E)$ be an undirected graph and $\lambda:V\to\mathbb{Z}$ a function
such that, for all $v\in V$,
\[0\le\lambda(v)\le\deg(v).\]

The following statements are equivalent:
\begin{itemize}
\item[(1)] For all subsets $A$ of $V$,
\[\sum_{v\in A}\lambda(v)\ge|E(G[A])|.\]
\item[(2)] There is an orientation $o$ of $G$ such that, for all $v\in V$,
\[\lambda(v)\ge\mathrm{In}^o(v),\]
the in-degree of $v$ in the orientation.
\end{itemize}

\paragraph{Question} What is the complexity of deciding whether these
conditions hold? (It is conjectured to be NP-hard.) In particular, can we do
better than $O(2^{|V|})$?

\subsection{Can you always beat a given cap set?}
\Ppresent{Fred Tyrrell}{fred.tyrrell@bristol.ac.uk}

A cap set is a subset $A$ of $\mathbb{F}_3^n$ such that $A$ has no non-trivial
solutions to $x+y+z=0$, other than when $x=y=z$. Equivalently, $A$ has no
non-trivial 3 term arithmetic progressions.

\paragraph{Problem} Given a cap set $A$ of size $c^n$ in $\mathbb{F}_3^n$, is
there some cap set $B$ in $\mathbb{F}_3^m$ for some $m>n$ such that
$|B| > c^m$?

\medskip

So far, the best known lower bounds for the cap set problem come from
constructing a fixed example, and then taking products of it with itself to
get asymptotic bounds. But it seems likely that the `true' lower bound does
not arise from a given cap set, but rather some family of cap sets, or perhaps
some limiting construction.

I can show the following sufficient condition for this question, using ideas
from my recent paper~\cite{tyrrell}.

Suppose $A$ and $A'$ are cap sets in $\mathbb{F}_3^n$, both of size $c^n$ for
some $c$. If $A$ and $A'$ are disjoint, then it is possible to construct a cap
set $B$ in $\mathbb{F}_3^m$ of size strictly greater than $c^m$, for some
$m>n$. In other words, if I can find two disjoint cap sets of the same size,
I can use these to construct a better cap set in some larger dimension.

However, I am not able to show that one can always find two disjoint cap sets
of the same size (and indeed it may well not be true!).

\subsection{The Four-Colour Theorem}
\Ppresent{Weiguo Xie}{xiew@d.umn.edu}

The four-color problem was first posed by Francis Guthrie in 1852. It remained
unsolved for over a century. Although computer-assisted proofs emerged in the 
last 50 years, they were typically the ``machine-checkable proofs'', which
were not therefore checked by human readers.

Is there an algorithm which, on input of a planar graph, outputs either a
$4$-coloring of it or ``fail'', with the properties
\clearpage
\begin{itemize}
\item it runs in polynomial time, and
\item the proportion of graphs on which it fails tends to $0$ as the number $n$
of vertices tends to infinity?
\end{itemize}
Even better would be an algorthm which never fails!

I have a procedure which may have this property. If anyone is interested,
please contact me (the presenter).

\end{document}